\documentclass[11pt]{article}

\usepackage{amssymb}
\usepackage{amsmath}
\usepackage{amsthm}

\usepackage{color}
\begin{document}

\parskip 1ex            
\parindent 5ex		

\renewcommand{\baselinestretch}{1.2}

\newtheorem{theorem}{Theorem}
\newtheorem{note}[theorem]{Note}
\newtheorem{lemma}[theorem]{Lemma}
\newtheorem{conjecture}[theorem]{Conjecture}
\newtheorem{proposition}[theorem]{Proposition}
\newtheorem{corollary}[theorem]{Corollary}
\newtheorem{claim}[theorem]{Claim}
\newtheorem{property}{Property}
\newtheorem{definition}{Definition}

\newcommand\mc{\mathcal}

\def\a{\alpha} 
\def\b{\beta} 
\def\d{\delta} 
\def\D{\Delta}
\def\e{\epsilon} 
\def\f{\phi} 
\def\F{{\Phi}} 
\def\vp{\varphi} 
\def\g{\gamma}
\def\G{\Gamma} 
\def\i{\iota} 
\def\k{\e/2} 
\def\K{\Kappa}
\def\z{\zeta} 
\def\th{\theta} 
\def\Th{\Theta}  
\def\l{\lambda}
\def\L{\Lambda} 
\def\m{\mu} 
\def\n{\nu} 
\def\p{\pi}
\def\r{\rho} 
\def\R{\Rho} 
\def\s{\sigma} 
\def\S{\Sigma}
\def\t{\tau} 
\def\om{\omega} 
\def\OM{\Omega} 
\def\Om{\Omega}
\def\U{\Upsilon} 
\def\l{\ell}

\def\P{\mathbb{P}}
\def\E{\mathbb{E}}

\newcommand\mb{\mathbb}

\parindent=0pt
\parskip=10pt

\title{Independence
 number of graphs with a prescribed number of cliques}
\author{
Tom Bohman\thanks{Department of Mathematical Sciences, Carnegie Mellon
University, Pittsburgh, PA 15213, USA. Email: {\tt tbohman@math.cmu.edu}.
Research supported in part by NSF grant DMS-1362785.}
\and
Dhruv Mubayi\thanks{Department of Mathematics, Statistics, and Computer Science, University of Illinois, Chicago, IL, 60607 USA.  Research partially supported by NSF grant DMS-1300138. Email: {\tt mubayi@uic.edu}}
}
\date{}
\maketitle

\begin{abstract} We consider the following problem posed by Erd\H os in 1962. Suppose that $G$ is an $n$-vertex  graph  where the number of $s$-cliques in $G$ is $t$. How small can the independence number of $G$ be?  Our main result suggests  that for fixed $s$, the smallest possible independence number  undergoes a transition at  $t=n^{s/2+o(1)}$. 

In the case of triangles ($s=3$) our method yields the following result which is sharp apart from constant factors and generalizes basic results in Ramsey theory: there exists $c>0$ such that every $n$-vertex graph with $t$ triangles has independence number at least
$$c \cdot \min\left\{ \sqrt {n \log n}\,  , \, \frac{n}{t^{1/3}} \left(\log \frac{n}{ t^{1/3}}\right)^{2/3} \right\}.$$
\end{abstract}

\section{Introduction}
An old problem in extremal graph theory due to Erd\H os~\cite{E62} is to determine $f(n,s,l)$, the minimum number of $s$-cliques over all graphs on $n$ vertices and independence number less than $l$. 
This problem has been extensively studied and has led to many other interesting questions. However, the range of parameters that have been investigated are confined to the case of fixed $l$~(see, e.g., \cite{DHMNS, PV, N}). 

 Here we consider the problem when  $s$ is fixed but $l$ can be a  power of $n$. It is more convenient to phrase the problem in inverse form as follows. Suppose that $s \ge 2$ is fixed and $G$ is an $n$-vertex graph in which the number of copies of $K_s$ is $t$, where $0 \le t \le {n \choose s}$. How small can the independence number of $G$ be? 
 
 When $s=2$ this question is answered by Tur\'an's theorem. However, when $s \ge 3$, one cannot hope to obtain a precise answer for this minimum for the full range of  $t$, since when $t=0$ it is equivalent to determining the Ramsey numbers $r(s,l)$ for fixed $s$ and large $l$. Indeed, our motivation for studying this parameter was to gain a better understanding of the Ramsey number $r(4,l)$. In trying to speculate about the growth rate of $r(4,l)$ a natural problem that arises is to determine the smallest possible independence number of a graph where we have a specified number of triangles (the case $s=3$ of our problem).
 
  Our first result below suggests that the behavior of this function changes at $t=n^{s/2+o(1)}$.
 
\begin{theorem} \label{main}
Let $s \ge 2$ be a fixed constant.  If $G$ is a graph on $n$ vertices 
with $t$ copies of $K_s$ then
\[ \alpha(G) \ge \begin{cases}
c_s n^{\frac{1}{s-1}} & \text{ if } t \le n^{s/2} \\
c_s \left( \frac{ n^s}{t} \right)^{\frac{1}{ \binom{s}{2} }} & \text{ if } t \ge n^{s/2} 
\end{cases} \]
\end{theorem}

\noindent
{\bf Remarks.}

\begin{enumerate}
\item The $s=2$ case of  Theorem~\ref{main} is (apart from constants) simply the bound given by  Tur\'an's Theorem.

\item The bound in the case of large $t$ is essentially sharp from random 
considerations. Indeed,  if $G$ is approximately random and has $t$ copies of
$K_s$ (where $t$ is sufficiently large) then
$G$ is roughly $ G_{n,p}$ with $ n^sp^{\binom{s}{2}} = \Theta( t ) $.  Then we expect
the independence number of $G$, up to a logarithmic factor, to be $1/p$.  This is exactly
the lower bound that we give in the large $t$ case.

\item Although our proof gives some additional $\log$ factors, we have chosen to omit those for clarity of presentation. For $s=3$ we  carry out the calculations in our method in further detail and provide more precise results below. 
\end{enumerate}
When $s=3$ Theorem~\ref{main} becomes
\[ \alpha(G) \ge \begin{cases}
c n^{\frac{1}{2}} & \text{ if } t \le n^{3/2} \\
c \frac{ n}{t^{1/3}}  & \text{ if } t \ge n^{3/2} 
\end{cases} \]
Here the bounds are sharp (apart from logarithmic factors) in the entire range of $t$ due to  the existence of Ramsey graphs with the desired independence number~\cite{r3t, BK2, FGM, K} for $t<n^{3/2}$ and approximately random graphs for $t>n^{3/2}$. However, our next result sharpens both upper and lower bounds considerably to obtain sharp results in order of magnitude.

\begin{theorem} \label{s=3}
There exists $c>0$ and $n_0$ such that  if $G$ is a graph on $n>n_0$ vertices 
with $t$ triangles  then
\[ \alpha(G) \ge \begin{cases}
c  \sqrt{n \log n}& \text{ if } t \le n^{3/2}\sqrt{\log n}  \\
c \frac{n}{t^{1/3}} \left(\log \frac{n}{t^{1/3}}\right)^{2/3} & \text{ if } t \ge n^{3/2}\sqrt{\log n}.
\end{cases} 
\]
Both bounds are sharp apart from the constant $c$. 
\end{theorem}

The related question of bounding the chromatic number of 
a graph with a fixed number of triangles was recently addressed by Harris \cite{Harris}.  We note
in passing that the two problems have different characters when the number of triangles is 
large.  The extremal graph for the chromatic number question in this regime is simply a clique while the
extremal graph for the independence number question is a direct product of a clique and a Ramsey graph.  Harris' proof from~\cite{Harris} uses a difficult theorem of Johansson as its main tool and in fact his bound on chromatic number together with a simple argument can be used to give another proof of Theorem~\ref{s=3}. Our proof uses only the lower bounds on independence number of triangle-free graphs; these are much easier results than the corresponding upper bounds on chromatic number.

\section{Proof of Theorem~\ref{main}}

 In the proof,  we carry the constant $ c_s$ through the calculation in the interest
of clarity.  All conditions that the constants must satisfy are listed as numbered
equations below.  We make no attempt to optimize these constants. 

We go by induction on $s$.  Let $d$ be the average degree of $G$.  
For an arbitrary vertex $x$ let $ d(x) $ denote
the degree of $x$ and let $ t(x) $ denote the number of copies of
$ K_s$ that contain $x$.

{\bf Case 1.} $ t \le n^{s/2} $.  

Here we show that $ \alpha(G) > c_s' n^{\frac{1}{s-1}} $ where $ c_s'$ is a constant (which will be
slightly larger than $c_s$ itself).

If $ d \le (c_s')^{-1} n^{ \frac{s-2}{s-1}}$ then we have the desired
bound by Tur\'an's Theorem:
\[ \alpha(G) \ge \frac{ n}{ d} \ge c_s' \frac{n}{ n^{ \frac{s-2}{s-1}} }
= c_s' n^{\frac{1}{s-1}}. \]
So, we may assume  $ d >  (c_s')^{-1} n^{ \frac{s-2}{s-1}} $. Now consider the random variable
$$ X_v =  d(v) - t(v)^{\frac{2}{s-1}} $$ where $v$ is a vertex chosen
uniformly at random.  We have
\begin{multline*}
\E\left[X_v\right] = 
\E \left[ d(v) - t(v)^{\frac{2}{s-1}} \right] \ge  (c_s')^{-1} n^{ \frac{s-2}{s-1}} - \E\left[ t(v) \right]^{\frac{2}{ s-1 }} \\
\ge  (c_s')^{-1} n^{ \frac{s-2}{s-1}} - \left( s n^{ \frac{s}{2} -1 } \right)^{\frac{2}{s-1}}
= \left( (c_s')^{-1}  - s^{\frac{2}{s-1}} \right)  n^{ \frac{ s-2}{s-1}} > 0.
\end{multline*}
Note that we assume
\begin{equation}
\label{eq:constant1}
c_s' < s^{- \frac{2}{s-1}}.
\end{equation}
Now we apply the inductive hypothesis in $G[N(v)]$ where $v$ is a vertex
that achieves the above bound on $X_v$.  This graph has 
$$ d(v) >  \left( (c_s')^{-1}  - s^{\frac{2}{s-1}} \right)  n^{ \frac{ s-2}{s-1}} $$
vertices and at most
$$ t(v) < d(v)^{ \frac{s-1}{2}}  $$
copies of $ K_{s-1}$.  Therefore
\[ \alpha(G) \ge \alpha( G[ N(v)] ) \ge c_{s-1} d(v)^{\frac{1}{s-2}} >  c_{s-1} \left( (c_s')^{-1}  -  s^{\frac{2}{s-1}} \right)^{\frac{1}{s-2}} n^{ \frac{1}{s-1}}. \] 
This gives the desired bound so long as
\begin{equation}
\label{eq:constant2}
c_s' \le  c_{s-1} \left( (c_s')^{-1}  -s^{\frac{2}{s-1}}  \right)^{\frac{1}{s-2}}. 
\end{equation}

{\bf Case 2.} $ t > n^{s/2}$.

Here we randomly sparsify our graph and apply Case 1.  We may assume $ t < \delta n^s $ where
$ \delta > 0 $ is a function of the constant $ c_s$ in the statement of the Theorem.

Set
\[ p = \frac{n}{ t^{2/s} 2^{1+ 2/s}} \]
and consider a set $S$ of vertices chosen at random so that each vertex appears
independently with probability $p$.  By standard concentration results (using the upper
$ t <  \delta n^s$ ) we 
have 
\[ |S| > \frac{np}{2} =  \frac{n^2}{ 2^{2+ 2/s}t^{2/s} }\]
with high probability.  
Let $T$ be the number of copies of $K_s$ in $ G[S]$.
By Markov's inequality, we have
\[ T \le 2 E[T] = 2 t p^s = \frac{ n^s  }{2^{s+1} t }  \] 
with probability at least $1/2$.  Thus, there exists a vertex set $S$ that satisfies
both of these inequalities.  So we can apply Case 1 to conclude
\[ \alpha (G) \ge \alpha ( G[S]) \ge  c_s' |S|^{\frac{1}{s-1}} \ge c_s' \frac{n^\frac{2}{s-1}}{ 2^{ \frac{2(s+1)}{s(s-1)}}t^{\frac{2}{s (s-1)}}}. \]
This gives the desired bound so long as
\begin{equation}
\label{eq:last}
c_s < c_s' 2^{ - \frac{2(s+1)}{s(s-1)}}.
\end{equation}

\section{Proof of Theorem~\ref{s=3}}
 We begin with the lower bounds. The following result will be used. It follows from the original proof of Ajtai-Koml\'os-Szemer\'edi~\cite{AKS} (see also Lemma 12.16 in~\cite{B}).
 
 \begin{theorem} \label{aks}
There exists an absolute constant $c$ such that the following holds. Let $G$ be an
$n$-vertex graph with average degree $d$ and $t$ triangles. Then
$$\alpha(G) > \frac{cn}{d}\left(\log d - \frac{1}{2} \log \left(\frac{t}{n}\right)\right).$$
\end{theorem}
\parindent=0pt

Now we continue with the lower bound proof of Theorem~\ref{s=3}. 

{\bf Case 1.} $t \le n^{3/2} (\log n)^{1/2}$.

Let $\epsilon=1/10$. If $d \le  n^{1/4 +\epsilon}$, then using Tur\'ans theorem we immediately obtain $\alpha(G) \ge n/(d+1)> c\sqrt{n \log n}$ with plenty of room to spare, so assume that $d>n^{1/4+\epsilon}$.  

Next suppose that $d>7\sqrt{n \log n}$.  Let $v$ be a vertex chosen uniformly at random.  Let $ d(v)$ be the degree of $v$  and let $t(v)$ be the number of 
triangles containing $v$. Then the expected value of $d(v)-2t(v)$ is at least $d-6t/n\ge \sqrt{n \log n}$ and since $\alpha(G) \ge d(v)-2t(v)$ we are done. Henceforth we may assume that $n^{1/4+\epsilon} < d < 7\sqrt{n \log n}$. Now we apply Theorem~\ref{aks} to obtain
$$\alpha(G) \ge \frac{cn}{d} \log\left(\frac{d\sqrt n}{\sqrt t}\right)> c' \sqrt{n \log n}$$
where $c'>0$.

{\bf Case 2.} $t > n^{3/2} (\log n)^{1/2}$.

Note that we can assume $ t < \delta n^3$ where $ \delta >0 $ is a function of the constant $c$ in the statement of the Theorem.  
To be precise, $ \delta$ can be taken to be solution of
the equation
$$ \frac{ c}{ \delta^{1/3}} \left( \log \frac{1}{ \delta^{1/3}} \right)^{2/3} = 1. $$

Set
$$p= \frac{1}{4} \left(  \frac{n}{t^{2/3}} \right) \left(\log \frac{n}{t^{1/3}} \right)^{1/3} <1.$$ Choose a vertex subset $S$ of $V(G)$ by picking each vertex independently with probability $p$. Then, by standard large deviation estimates for $np$ large and standard Poisson approximation 
estimates for $np$ constant, with high probability $$|S| \ge pn/2 :=N.$$  And by Markov's inequality, with probability at least 1/2, the number of triangles $T$ in $S$ satisfies
$$T<2p^3t= \frac{1}{32} \left( \frac{n^3}{t} \right) \left(\log \frac{n}{ t^{1/3}} \right)^{1/3}< N^{3/2}\sqrt{\log N}.$$
Thus, we can choose a set $S$ that satisfies both of these inequalities and 
apply the result we have already proved in Case 1 to obtain the lower bound
\begin{multline*}
\alpha(G) \ge \alpha(G[S]) \ge c'\sqrt{N \log N} \ge c'\sqrt{\frac{n^2(\log (n/ t^{1/3}) )^{1/3}  (2 \log (n/t^{1/3}))}{8 t^{2/3}}} \\ = \frac{c'}{2} \left( \frac{n}{t^{1/3}} \right) \left(\log \frac{n}{t^{1/3}} \right)^{2/3}.
\end{multline*}
Note that the last inequality holds for $ \delta$, and hence $c$ sufficiently small.

Now we will exhibit constructions showing that these bounds are sharp. In the range  $t \le n^{3/2} (\log n)^{1/2}$ we proceed as follows. Choose $a$ so that ${a \choose 3} =t$ and let $G'$ be a triangle-free graph on $n-a$ vertices with independence number at most $c\sqrt{(n-a)\log(n-a)}$ (this exists due to~\cite{K}). Note that the upper bound on $t$ implies that $a<2\sqrt n \log n= o(n)$.
Let $G$ be the graph that is a disjoint union of $G'$ and $K_a$. Then $G$ has $n$ vertices, $t$ triangles, and $\alpha(G) =\alpha(G')+1$.

Next we consider the range $t \ge n^{3/2} (\log n)^{1/2}$.  We may assume $ t < \delta n^3$ for some $ \delta >0 $. 
Define $\lambda>0$ by the equation
$$t=n^{3/2} \lambda^{3/2} \sqrt{\log \frac{n}{ t^{1/3}}}$$ and assume for simplicity of notation that $\lambda$ is an integer. Let $N=n/\lambda$ and again assume that $N$ is an integer. Let $G'$ be a triangle-free graph on $N$ vertices and independence number at most $c\sqrt{N \log N}$. Note that $G'$ has $\Theta(N^{3/2}\sqrt{\log N})$ edges. Let $G$ be the lexicographic product of $G'$ with $K_{\lambda}$. In other words, replace each  vertex $v$ of $G'$ with a clique $S_v$ of size $\lambda$ and for $x \in S_v$, we have $N_G(x) = (S_v\setminus \{v\})\cup \cup_{w \in N_{G'}(v)} S_w$.  Then $G$ has $\lambda N=n$ vertices and the number of triangles in $G$ is
$t'=\Theta(\lambda^3|E(G')|) = \Theta(t)$. Finally, 
$$\alpha(G)=\alpha(G')\le c\sqrt{N \log N}  < c' \frac{n}{t^{1/3}} (\log (n/ t^{1/3}))^{2/3}.$$
We have produced a construction where the number of triangles is $t'=\Theta(t)$ and it is an easy matter to produce one where the number of triangles is exactly $t$ and the independence number has the same order of magnitude. \qed 

{\bf Acknowledgment.}   The authors thank David Harris for pointing out that the results in~\cite{Harris} can be adapted to yield another proof of Theorem~\ref{s=3}. The The first author thanks Gregory and Nina Goossens.


\begin{thebibliography}{99}

\bibitem{AKS} M. Ajtai, J. Koml\'os, and E. Szemer\'edi. A note on Ramsey numbers. J. Combin. Theory
Ser. A, 29(3): 1980 354–-360.


\bibitem{r3t} T. Bohman, The triangle-free process, {\em Advances in Mathematics}
{\bf 221} (2009) 1653--1677.


\bibitem{BK2} T. Bohman and P. Keevash,
Dynamic concentration of the triangle-free process.
{\em arXiv:1302.5963}

\bibitem{B} B. Bollob\'as. Random graphs, volume 73 of Cambridge Studies in Advanced Mathematics.
Cambridge University Press, Cambridge, second edition, 2001.

\bibitem{DHMNS} S. Das,  H.  Huang,  J. Ma, H. Naves, B.  Sudakov, 
A problem of Erdős on the minimum number of k-cliques. (English summary) 
J. Combin. Theory Ser. B 103 (2013), no. 3, 344–-373. 


\bibitem{E62} P. Erd\H os,
On the number of complete subgraphs contained in certain graphs,
Magyar Tud. Akad. Mat. Kutató Int. Közl. 7 1962 459-–464. 


\bibitem{FGM} G. Fiz Pontiveros, S. Griffiths and R. Morris,  
The triangle-free process and R(3,k). {\em arXiv:1302.6279}

\bibitem{Harris} D. Harris,
Some results on chromatic number as a function of triangle count.
{\em arXiv:1604.00438}

\bibitem{K} J.H. Kim,
The Ramsey number $R(3,t)$ has order of magnitude $t^2/\log t$,
{\em Random Structures Algorithms} {\bf 7} (1995), 173--207.

\bibitem{N} V. Nikiforov,  On the minimum number of k-cliques in graphs with restricted independence number. Combin. Probab. Comput. 10 (2001), no. 4, 361-–366. 

\bibitem{PV} O. Pikhurko, E.  Vaughan,
Minimum number of k-cliques in graphs with bounded independence number, 
Combin. Probab. Comput. 22 (2013), no. 6, 910–-934. 


\end{thebibliography}
\end{document}